\documentclass[12pt,draft]{amsart}
\usepackage{amsmath,amsthm,latexsym,amscd,amsbsy,amssymb,fleqn,leqno}
\chardef\bslash=`\\ 

\makeatletter
\def\verbatim{\interlinepenalty\@M \@verbatim
  \leftskip\@totalleftmargin\advance\leftskip2pc
  \frenchspacing\@vobeyspaces \@xverbatim}
\makeatother
\makeatletter
  \def\dgt@k{\dg@DX=-3 \dg@DY=2 \dg@SIZE=3} 
\makeatother

\makeatletter
  \def\dgt@kk{\dg@DX=3 \dg@DY=-1 \dg@SIZE=3}%
\makeatother

\theoremstyle{plain}
\newtheorem{thm}{Theorem}[section]

\newtheorem{pro}[thm]{Proposition}

\theoremstyle{definition}

\numberwithin{equation}{section}

\newcounter{rmnum}

\def\symbolnote#1#2{\let\thefootn=\thefootnote%
\renewcommand{\thefootnote}{\fnsymbol{footnote}}%
\footnotemark[#1]%
\footnotetext[#1]{#2}%
\let\thefootnote=\thefootn
}

\newfont{\bbb}{msbm10 scaled \magstep1}
\newfont{\bbc}{msbm8 scaled \magstep0}

\newcommand{\R}{\mbox{\bbb R}}

\begin{document}


\title{On intersection of simply connected sets in the plane }

\author{E. D. Tymchatyn}
\address{Department of Mathematics and Statistics, University of Saskatchewan,
McLean Hall, 106 Wiggins Road, Saskatoon, SK, S7N 5E6, Canada}
\email{tymchat@math.usask.ca}
\thanks{The first author was partially supported by NSERC grant NPG 0005616}
\thanks{The second author was partially supported by NSERC grant RGPIN 261914}

\author{Vesko Valov}
\address{Department of Computer Science and Mathematics, Nipissing University,
100 College Drive, P.O. Box 5002, North Bay, ON, P1B 8L7, Canada}
\email{veskov@nipissingu.ca}

 

\begin{abstract}
Several authors \cite{bo} and \cite{kr} have recently attempted to show that the intersection of three simply connected subcontinua of the plane is simply connected provided it is non-empty and the intersection of each two of the continua is path connected. In this note we give a very short complete proof of this fact. We also confirm a related conjecture of Karimov and Repov\v{s} \cite{kr}. 
\end{abstract}

\maketitle

\markboth{E. Tymchatyn and V.~Valov}{On intersection of simply connected sets }


\section{Introduction}

A homology (resp., singular) cell is a compact metric space whose Vietoris (resp., singular) homology groups are trivial.
Helly \cite{he} proved the following result which is now known as the Topological Helly Theorem:

\begin{thm}
Let $\displaystyle {\mathcal S}=\{S_0,...,S_m\}$, $m\geq n$, be a finite family of homology cells in $\R^n$ such that the intersection of every subfamily ${\mathcal H}$ of ${\mathcal S}$ is nonempty if the cardinality $|\mathcal H|\leq n+1$ and it is a homology cell if $|\mathcal H|\leq n$. Then $\cap_{i=0}^{i=m}S_i$ is a homology cell.  
\end{thm} 

\medskip
Versions of Theorem 1.1 for singular homology have been proved by Debrunner \cite{d}  and  Alexandroff and   Hopf  \cite[p. 295]{ah}
for open sets in $\R^n$ and simplicial complexes in $\R^n$, respectively.

\medskip
A topological space is said to be  simply connected if it is path connected and has trivial fundamental group. 
 It is known  \cite{ccz} that a compact subspace of the plane is a singular cell if and only if it is simply connected. 

\medskip
In section 2 of the paper \cite{he} Helly proved that if $S_i$, $i=1,..,4$, are singular cells in $\R^2$ such that all intersections 
$\displaystyle S_{i_1}\cap S_{i_2}\cap S_{i_3}$ are singular cells, then $\cap_{i=1}^{i=4}S_i$ is not empty. Hence to prove the Topological Helly Theorem for singular cells in $\R^2$, it suffices to prove the following: 

\begin{pro}
Let $S_0, S_1$  and $S_2$ be three simply connected compacta in the plane such that the intersection of any two of them is path connected and  $\displaystyle \cap_{i=0}^{i=2} S_i\neq\emptyset$. 
Then $\displaystyle \cap_{i=0}^{i=2}S_i$ is simply connected.
\end{pro}

Bogatyi \cite{bo} has pointed out that no complete proof of this proposition can be found in the literature. He proved the proposition in the special case that $S_i$ are Peano continua. Karimov and Repov\v{s} \cite{kr}, established that, with the hypotheses of Proposition 1.2,     
$\displaystyle \cap_{i=0}^{i=2}S_i$ is cell-like connected (i.e., every two points can be connected by a cell-like continuum).
We prove Proposition 1.2  by showing that $\displaystyle \cap_{i=0}^{i=2}S_i$ is path connected. We also give an affirmative answer to a conjecture of Karimov and Repov\v{s} \cite{kr} by proving the following proposition:

\begin{pro}
If $X$ and $Y$ are compact $AR$'s in the plane, then so is each component of $X\cap Y$.
\end{pro} 

\medskip
We would like to thank the referee for his/her valuable remarks. We also thank Dr. U. Karimov for pointing out an error in the first version of this note and Dr. S. Bogatyi for several fruitful discussions and comments.


\section{Proof of Proposition 1.2}

Since the intersection of any family of simply connected sets in the plane has a trivial fundamental group with respect to each of its points, it suffices to show that 
$\displaystyle \bigcap_{i=0}^{i=2}S_i$ is path connected.
Suppose this is not the case and let 
$0$ and $1$ be two points in distinct arc components of $\displaystyle \bigcap_{i=0}^{i=2}S_i$. Let $\displaystyle I\subset S_0\cap S_1$,  
$\displaystyle J\subset S_0\cap S_2$ and $\displaystyle K\subset S_1\cap S_2$ be arcs from $0$ to $1$. Consider the components $J_n$, $n=1,2,..$, of $\displaystyle J\backslash\big(I\cup K\big)$ which are not in $S_1$. Since $0$ and 1 are end-points of $J$ it follows that no $J_i$ separates $\displaystyle I\cup J\cup K$. Also, no $J_i$ lies in a bounded component of $\displaystyle\R^2\backslash\big(I\cup K\big)$ 
because if the locally connected continuum $I\cup K$ separates $J_i$ from $\infty$ in $\R^2$,
then some simple closed curve in $\displaystyle I\cup K\subset S_1$ would do so as well. Since $S_1$ is simply  
connected, this would imply $\displaystyle J_i\subset S_1$. 

We are going to construct for every $n\geq 1$ an arc $\displaystyle J^n\subset S_0\cap S_1$ in $\displaystyle I\cup J\cup K$ from $0$ to 1 such that $\displaystyle J^n\cap\big(J_1\cup..\cup J_n\big)=\emptyset$. Let $D_1$ be the bounded component of  
$\displaystyle\R^2\backslash (I\cup J\cup K)$ whose boundary contains $J_1$ and 
$J_1\cup I\cup K$ separates $D_1$ from infinity in $\R^2$. Then  
$\displaystyle D_1\subset D_I\cap D_K$, where $D_I$ (resp., $D_K$) is the component of $\displaystyle\R^2 \backslash (I\cup J)$ (resp., $\displaystyle\R^2\backslash (J\cup K)$) containing $D_1$. Then 
$D_I$ and $D_K$ are bounded because $J_1$ does not separate $\displaystyle I\cup J\cup K$. Note that 
$\displaystyle D_I\subset S_0$ and $\displaystyle D_K\subset S_2$. Let $D$ be the component of $\displaystyle D_I\cap D_K$ containing $D_1$. Then $\displaystyle D\subset S_0\cap S_2$ and $\displaystyle J_1\subset\overline{D}\subset S_0\cap S_2$. Moreover, $\displaystyle Fr(D)\subset I\cup J\cup K$. It is well known \cite{St} that each continuum contained in the union of finitely many arcs is rim-finite and, hence, locally connected.  So $Fr(D)$ is locally connected. Let $C\subset Fr(D)$ be the simple closed curve that separates $D$ from $\infty$ in $\R^2$. Then there is an arc $\displaystyle J^1\subset\big(J\cup C\big)\backslash J_1\subset S_0\cap S_2$ from $0$ to $1$. Obviously, $\displaystyle J^1\subset\R^2\backslash J_1$ since $J_1\subset D$.
Suppose we already constructed an arc $\displaystyle J^n\subset S_0\cap S_1$ in $\displaystyle I\cup J\cup K$ from $0$ to 1 such that $\displaystyle J^n\cap\big(J_1\cup..\cup J_n\big)=\emptyset$. If $\displaystyle J^n\cap J_{n+1}=\emptyset$, let $J^{n+1}=J^n$. If  $\displaystyle J^n\cap J_{n+1}\neq\emptyset$, we repeat the above arguments with $\displaystyle J^n$ in place of $J$ and $\displaystyle J_{n+1}$ in place of $\displaystyle J_1$ to obtain an arc 
$\displaystyle J^{n+1}\subset S_0\cap S_2\cap\big(J^n\cup I\cup J\backslash J_{n+1}\big)$ from 0 to 1. By induction, we construct a sequence of arcs $\displaystyle\{J^n\}_{n=1}^{\infty}$ from 0 to 1 with 
$\displaystyle J^{n+1}\subset S_0\cap S_1\cap\big(I\cup J\cup K\backslash\bigcup_{i=1}^{n+1}J_i\big)$.   
Let $\displaystyle J^*=\limsup J^n$. Then $\displaystyle J^*\subset\big(S_0\cap S_2)\cap\big(I\cup J\cup K\backslash\bigcup_{i=1}^{\infty}J_i\big)\subset S_1$ is a continuum from 0 to 1. As above, $J^*$ is locally connected. So, there is an arc in $J^*$ from 0 to 1 which contradicts the fact that $0$ and $1$ are in distinct arc components of $\displaystyle \bigcap_{i=0}^{i=2}S_i$.


\section{Proof of Proposition 1.3}

Let $C$ be a component of $X\cap Y$. If $K$ is the topological hull of $C$, then $K\subset X$ and $K\subset Y$ since neither $X$ nor $Y$ separates $\R^2$. So, $K=C$. By unicoherence of $\R^2$ it follows that $Fr(C)$, the boundary of $C$ in $\R^2$, is connected. 

By the well-known result of Borsuk \cite{br} (that every locally connected plane continuum not separating the plane is an $AR$), it remains to prove that $C$ is locally connected.  Since $C$ is a continuum in the plane, it suffices to prove that  $Fr(C)$ is locally connected. To prove this it suffices to show that every pair of points of $Fr(C)$ is separated by a finite set (see \cite[p. 99]{wh}).

Since $X$ is simply connected, locally connected subcontinuum in the plane, by \cite[ ch. IV]{wh}, all true cyclic elements of $X$ are topological disks $D_i$ such that  the cardinality of $D_i\cap D_j$ is at most 1 for $i\neq j$  and, if the sequence $\{D_i\}$ is infinite, then $\lim diam D_i=0$  .  Hence, 
each $Fr(D_i)$ is a simple closed curve and $Fr(X)=X\backslash\bigcup int(D_i)$ is a locally connected continuum with a particularly simple structure. Let $x$ and $y$ be distinct points in $Fr(C)\subset Fr(X)\cup Fr(Y)$. If $x$ and $y$ do not both lie in any one cyclic element of $X$, then an one point set separates $x$ and $y$ in $X$ and, hence, in $C$. Thus, we may suppose that there are cyclic elements $D$ in $X$ and $E$ in $Y$ with $x,y\in D\cap E$.
Now $x$ in $int(D)$ implies there is a neighborhood $W$ of $x$ in $Fr(X)\cup Fr(Y)$ 
with $\overline{W}\subset int(D)$. Then a finite set $P$ separates $Fr(Y)\backslash W$ from $x$ in $Fr(Y)$ since $Fr(Y)$ is rim-finite. Hence, $P$ separates $x$ from $Fr(X)\cup Fr(Y)\backslash W$. So we may suppose $x,y\in Fr(D)\cap Fr(E)$ (see \cite[49.V, Theorem 3, p. 244]{k}).

Let $F$ be a two-point set in $Fr(E)$ which separates $x$ and $y$ in $Fr(E)$. Then $F$ separates $x$ and $y$ in $Fr(Y)$ \cite[IV.3.1, p. 67]{wh}. Since $D$ is hereditary normal, there is a closed set $A\subset D$ which separates $x$ and $y$ in $D$ and such that $A\cap Y\subset F$. Since $D$ is unicoherent, a component $A^{'}$ of $A$  separates $x$ and $y$ in $D$. It is now a routine exercise to construct an arc $A^{''}\subset D$ such that $A^{''}$ separates $x$ and $y$ in $D$ and $A^{''}\cap Y\subset F$. If we also take $A^{''}$ to be
irreducible with respect to separating $x$ and $y$ in $D$ (see \cite[V.49, Theorem 3, p.244]{k}), then $A^{''}\cap Fr(D)$ will contain just two points $c$ and $d$. As above, 
$A^{''}$ separates $x$ and $y$ in $X$ because $D$ is a cyclic element of $X$.
So $A^{''}\cap\big(Fr(X)\cup Fr(Y)\big)\subset F\cup\{c,d\}$ separates $x$ and $y$ in 
$Fr(C)\subset\big(Fr(X)\cup Fr(Y)\big)\subset X$ .So, $Fr(C)$ is rim-finite, hence, locally 
connected.

\end{document}